\documentclass[12pt]{article}

\usepackage{amsmath, epsfig, cite}
\usepackage{amssymb}
\usepackage{amsfonts}
\usepackage{latexsym}
\usepackage{graphicx}
\usepackage{color}
\usepackage{ifpdf}

\newtheorem{thm}{Theorem}[section]

\newtheorem{lem}[thm]{Lemma}
\newtheorem{conj}[thm]{Conjecture}

\def\pf{\noindent {\it Proof.} }

\numberwithin{equation}{section}

\makeatletter \@addtoreset{equation}{section} \makeatother

\begin{document}
\rule{0cm}{1cm}

\begin{center}
{\Large\bf The (revised) Szeged index and the\\[2mm] Wiener index of a
nonbipartite graph }
\end{center}

\begin{center}
{\small Lily Chen, Xueliang Li, Mengmeng Liu\\
Center for Combinatorics, LPMC-TJKLC\\
Nankai University, Tianjin 300071, China\\
Email:  lily60612@126.com, lxl@nankai.edu.cn, liumm05@163.com}
\end{center}

\begin{center}
\begin{minipage}{120mm}
\begin{center}
{\bf Abstract}
\end{center}

{\small Hansen et. al. used the computer programm AutoGraphiX to
study the differences between the Szeged index $Sz(G)$ and the
Wiener index $W(G)$, and between the revised Szeged index $Sz^*(G)$
and the Wiener index for a connected graph $G$. They conjectured
that for a connected nonbipartite graph $G$ with $n \geq 5$ vertices
and girth $g \geq 5,$ $ Sz(G)-W(G) \geq 2n-5. $ Moreover, the bound
is best possible as shown by the graph composed of a cycle on $5$
vertices, $C_5$, and a tree $T$ on $n-4$ vertices sharing a single
vertex. They also conjectured that for a connected nonbipartite
graph $G$ with $n \geq 4$ vertices, $ Sz^*(G)-W(G) \geq
\frac{n^2+4n-6}{4}. $ Moreover, the bound is best possible as shown
by the graph composed of a cycle on $3$ vertices, $C_3$, and a tree
$T$ on $n-3$ vertices sharing a single vertex. In this paper, we not
only give confirmative proofs to these two conjectures but also
characterize those graphs that achieve the two lower bounds. }

\vskip 3mm

\noindent {\bf Keywords:} Wiener index, Szeged index, revised Szeged
index. \vskip 3mm


\end{minipage}
\end{center}

\section{Introduction}

All graphs considered in this paper are finite, undirected and
simple. We refer the readers to \cite{bm} for terminology and
notation. Let $G$ be a connected graph with vertex set $V(G)$ and
edge set $E(G)$. For $u,v \in V$, $d_G(u,v)$ denotes the {\it
distance} between $u$ and $v$ min $G$. The {\it Wiener index} of $G$
is defined as
$$
W(G)=\displaystyle\sum_{\{u,v\}\subseteq V(G)} d_G(u,v).
$$
This topological index has been extensively studied in the
mathematical literature; see, e.g., \cite{GSM,GYLL}. Let $e=uv$ be
an edge of $G$, and define three sets as follows:
$$
N_u(e) = \{w \in V: d_G(u,w)< d_G(v,w)\},
$$
$$
N_v(e) = \{w \in V: d_G(v,w)< d_G(u,w)\},
$$
$$
N_0(e) = \{w \in V: d_G(u,w)=d_G(v,w)\}.
$$
Thus, $\{N_u(e),N_v(e),N_0(e)\}$ is a partition of the vertices of
$G$ respect to $e$. The number of vertices of $N_u(e)$, $N_v(e)$ and
$N_0(e)$ are denoted by $n_u(e)$, $n_v(e)$ and $n_0(e)$,
respectively. A long time known property of the Wiener index is the
formula \cite{GP,W}:
$$
W(G) = \displaystyle\sum_{e=uv \in E(G)} n_u(e) n_v(e),
$$
which is applicable for trees. Motivated the above formula, Gutman
\cite{G} introduced a graph invariant, named as the {\it Szeged
index} as an extension of the Wiener index and defined by
$$
Sz(G) = \displaystyle\sum_{e=uv \in E(G)}n_u(e) n_v(e).
$$
Randi\'c \cite{R} observed that the Szeged index does not take into
account the contributions of the vertices at equal distances from
the endpoints of an edge, and so he conceived a modified version of
the Szeged index which is named as the {\it revised Szeged index}.
The revised Szeged index of a connected graph $G$ is defined as
$$
Sz^*(G) = \displaystyle\sum_{e=uv \in E(G)}\left(n_u(e)+
\frac{n_0(e)}{2}\right)\left(n_v(e)+ \frac{n_0(e)}{2}\right).
$$

Some properties and applications of the Szeged index and the revised
Szeged index have been reported in \cite{AH,CLL,LL,PR,PZ,XZ}.

We have known that for a connected graph, $Sz^*(G) \geq Sz(G) \geq
W(G)$. It is easy to see that $Sz^*(G)=Sz(G)=W(G)$ if $G$ is a tree,
which means $m=n-1$.  So we want to know the differences between
$Sz(G)$ and $W(G)$, and between $Sz^*(G)$ and $W(G)$ for a connected
graph with $m \geq n$.

In \cite{auto} Hansen et. al. used the computer programm AutoGraphiX
and made the following conjectures:

\begin{conj}\label{conj1}
Let $G$ be a connected bipartite graph with $n \geq 4$ vertices and
$m \geq n$ edges. Then
$$
Sz(G)-W(G) \geq 4n-8.
$$
Moreover, the bound is best possible as shown by the graph composed
of a cycle on $4$ vertices $C_4$ and a tree $T$ on $n-3$ vertices
sharing a single vertex.
\end{conj}

\begin{conj}\label{conj2}
Let $G$ be a connected bipartite graph with $n \geq 4$ vertices and
$m \geq n$ edges. Then
$$
Sz^*(G)-W(G) \geq 4n-8.
$$
Moreover, the bound is best possible as shown by the graph composed
of a cycle on $4$ vertices $C_4$ and a tree $T$ on $n-3$ vertices
sharing a single vertex.
\end{conj}

\begin{conj}\label{conj3}
Let $G$ be a connected graph with $n \geq 5$ vertices and girth
$g\geq 5$ and with an odd cycle. Then
$$
Sz(G)-W(G) \geq 2n-5.
$$
Moreover, the bound is best possible as shown by the graph composed
of a cycle on $5$ vertices $C_5$ and a tree $T$ on $n-4$ vertices
sharing a single vertex.
\end{conj}

\begin{conj}\label{conj4}
Let $G$ be a connected graph with $n \geq 4$ vertices and $m \geq n$
edges and with an odd cycle. Then
$$
Sz^*(G)-W(G) \geq \frac{n^2+4n-6}{4}.
$$
Moreover, the bound is best possible as shown by the graph composed
of a cycle on $3$ vertices $C_3$ and a tree $T$ on $n-3$ vertices
sharing a single vertex.
\end{conj}

In \cite{CLL} we showed that both Conjecture \ref{conj1} and
\ref{conj2} are true. In this paper, we will give confirmative
proofs to Conjecture \ref{conj3} and Conjecture \ref{conj4}. During
the proof of Conjecture \ref{conj3}, we find another case which also
makes the equality holds, that is the graph composed of a cycle on
$5$ vertices, $C_5$, and two trees with roots $v_1,v_2$ in $C_5$
satisfying $v_1v_2 \in E(C_5)$. So we get the following theorem:
\begin{thm}\label{thm1}
Let $G$ be a connected nonbipartite graph on $n \geq 5$ vertices and
girth $g\geq 5$. Then
$$
Sz(G)-W(G) \geq 2n-5.
$$
Equality holds if and only if $G$ is composed of a cycle $C_5$ on
$5$ vertices, and one tree rooted at a vertex of the cycle $C_5$ or
two trees, respectively, rooted at two adjacent vertices of the
cycle $C_5$.
\end{thm}

We notice that the method used in the proof of Theorem \ref{thm1}
can also be used to prove the bipartite case, and therefore this
gives another proof to Conjecture \ref{conj1} other than that in
\cite{CLL}.

\section{Main results}

We start this section with two definitions that are needed in our
later proofs frequently.

\textbf{Definition 1.} Let $P$ be a shortest path between two
vertices $x$ and $y$ in a graph $G$, $P'$ another path from $x$ to
$y$ in $G$. We call $P'$ the \emph{second shortest path} between $x$
and $y$, if $P'\neq P$, $|P'|-|P|$ is minimum, and if there are more
than one path satisfying the condition, we choose $P'$ as a one with
the most common vertices with $P$ in $G$.

\textbf{Definition 2.} A subgraph $H$ of a graph $G$ is called
\emph{isometric} if distance between any pair of vertices in $H$ is
the same as their distance in $G.$

In \cite{SGB} Gutman gave another expression for the Szeged index:
$$
Sz(G)= \displaystyle\sum_{e=uv \in E(G)}n_u(e) n_v(e)=
\displaystyle\sum_{e=uv \in E(G)}\displaystyle\sum_{\{x,y\} \subseteq
V(G)}\mu_{x,y}(e)
$$
where $\mu_{x,y}(e)$, interpreted as contribution of the vertex pair
$x$ and $y$ to the product $n_u(e)n_v(e)$, is defined as follows:

$$
\mu_{x,y}(e) =  \left\{
\begin{array}{ll}
1,& \mbox {if}
\begin{cases} d_G(x,u)<d_G(x,v) \text{ and } d_G(y,v)<d_G(y,u), \\ \text{or} \\
d_G(x,v)<d_G(x,u)\text{ and }\ d_G(y,u)<d_G(y,v),\\
\end{cases}\\
0, & \mbox{otherwise}.
\end{array}
\right.
$$

From above expressions, we know that
\begin{eqnarray*}
Sz(G)-W(G) & = & \displaystyle \sum_{\{x,y\} \subseteq V(G)}
\displaystyle \sum_{e \in E(G)}\mu_{x,y}(e)-
\displaystyle \sum_{\{x,y\} \subseteq V(G)}d_G(x,y)\\
&=& \displaystyle \sum_{\{x,y\} \subseteq V(G)}\left( \displaystyle
\sum_{e \in E(G)}\mu_{x,y}(e)-d_G(x,y)\right).
\end{eqnarray*}

For convenience, let $\pi(x,y)=\sum_{e \in E(G)}\mu_{x,y}(e)-d_G(x,y)$.

Let $G$ be a connected graph. For every pair $\{x,y\} \subseteq
V(G)$, let $P_1$ be the shortest path between $x$ and $y$. We know
that for all $e \in E(P_1), \mu_{x,y}(e)=1$, which means that
$\pi(x,y)\geq 0.$ Let $P_2$ be the second shortest path between $x$
and $y$ (if there exists). Then $P_1 \Delta P_2 =C,$ where $C$ is a
cycle. Let $P'_i=P_i\bigcap C=x'P_iy'$. If $E(P_1)\bigcap
E(P_2)=\emptyset,$ then $x'=x, y'=y.$

Now we have the following lemma.

\begin{lem}\label{lem1}
For every pair $\{x,y\} \subseteq V(G),$ and $C,x',y'$ defined as above,\\
(1) if $C$ is an even cycle, then $\pi(x,y)\geq d_C(x',y')\geq 1;$\\
(2) if $C$ is an odd cycle and $d_C(x',y')\geq 2,$ then
$\pi(x,y)\geq 1.$
\end{lem}

\pf Firstly, we prove that for every $v \in V(C),$
$d_C(x',v)=d_G(x',v).$ If $v \in P'_1$, it is simply true;
otherwise, we can find a shorter path between $x'$ and $y'$, and
then we can find a shorter path between $x$ and $y$. If $v \in P'_2$
and $d_C(x',v) > d_G(x',v) =|E(P_3)|$, where $P_3$ is a shortest
path between $x'$ and $v$ in $G$, then the path
$xP_2x'P_3vP_2y'P_2y$ between $x$ and $y$ is shorter than $P_2$, a
contradiction. For the same reason, we have  $d_C(y',v)=d_G(y',v)$
for all $v \in V(C).$ Similarly, it is easy to see that a shortest
path from $x$ (or $y$) to the vertex $v$ in $P'_2$ is
$xP_2x'(yP_2y')$ together with a shortest path from $x'(y')$ to $v$
in $C$. So, if an edge $e$ in $E(C)$ makes $\mu_{x',y'}(e) =1$, then
we have $\mu_{x,y}(e) =1$.

If $C$ is an even cycle, we know that $|E(P'_2)|\geq |E(P'_1)|.$ For
every edge $e$ in the antipodal edges of $P'_1$ in $C$, it is
obviously that $\mu_{x',y'}(e) =1,$ and then $\mu_{x,y}(e) =1.$
Hence, $\sum_{e \in E(G)}\mu_{x,y}(e)=d_G(x,y)+ d_C(x',y')$, which
means that $\pi(x,y)\geq d_C(x',y')\geq 1.$

If $C$ is an odd cycle, there are vertices $x_1,x_2,y_1,y_2$ such
that
\begin{eqnarray}\label{equ1}
d_C(x',x_1)=d_C(x',x_2),\\ \label{equ2}
d_C(y',y_1)=d_C(y',y_2).
\end{eqnarray}
Let $d_C(x_1,y_1)=min \{\ d_C(x_i,y_j),i,j\in\{1,2\}\}$.  For every
edge $e$ in a shortest path between $x_1$ and $y_1$, we have
$\mu_{x,y}(e) =\mu_{x',y'}(e)=1$. So,  $\sum_{e \in
E(G)}\mu_{x,y}(e)\geq d_G(x,y)+ d_C(x_1,y_1)$, which means that
$\pi(x,y)\geq d_C(x_1,y_1).$

Next we show that $d_C(x_1,y_1)\geq 1$. From equations \ref{equ1}
and \ref{equ2}, we have
$$
d_C(x',x_1)=d_C(x',y')+d_C(y',x_1)-1,
$$
$$
d_C(y',y_1)=d_C(x',y')+d_C(x',y_1)-1.
$$
If $d_C(x_1,y_1)=0$, that is $x_1=y_1$, then by adding the above two
equations, we get
$$
d_C(x',y')=1,
$$
which contradicts the assumption $d_C(x',y')\geq 2.$

\begin{qed}
\end{qed}

From the proof of Lemma \ref{lem1}, we also get the following lemma.

\begin{lem}\label{lem2}
For every pair $\{x,y\} \subseteq V(C),$ where $C$ is an isometric cycle,\\
(1) if $C$ is an even cycle, then $\pi(x,y)\geq d_C(x,y)\geq 1;$\\
(2) if $C$ is an odd cycle and $d_C(x,y)\geq 2,$ then $\pi(x,y)\geq
1.$
\end{lem}

Now, we give a confirmative proof of Theorem \ref{thm1}.

\textbf{Proof of Theorem \ref{thm1}:} Let $C=v_1v_2\cdots v_kv_1$ be
a shortest odd cycle of $G$ with length $k$, where $k\geq g\geq 5$.
It is obvious that $C$ is an isometric cycle. We consider the pair
$\{x,y\}\subseteq V(G)$.

\textbf{Case 1.}  $\{x,y\}\subseteq V(C)$.

If $d_C(x,y) \geq 2$, then by Lemma \ref{lem2} we have $\pi(x,y)\geq
1$. Otherwise, $\pi(x,y)\geq 0$. Therefore,
$$\sum_{\{x,y\}\subseteq V(C)}\pi(x,y)\geq {k\choose 2} -k.$$

\textbf{Case 2.}  $x\in V(C), y\in V(G)\backslash V(C)$.

We will prove that for every $y \in V(G)\backslash V(C)$,
there exist two vertices $x_1, x_2$ in $C$ such that
$\pi(x_1,y)\geq 1$ and $\pi(x_2,y)\geq 1$.

Assume that $v_i$ is the vertex on $C$ such that
$d_G(v_i,y)=min_{v\in V(C)}d_G(v,y)$, and $P_1$ is a shortest path
between $v_i$ and $y$. Let $|E(P_1)|=p_1$. It is obvious that $P_1$
does not contain any vertex in $C$.

Now we show that $\pi(v_{i+2},y)\geq 1$. Since $P_2=yP_1v_iv_{i+1}v_{i+2}$ is a path
from $y$ to $v_{i+2}$, $p_1=d_G(v_i,y)\leq d_G(v_{i+2},y)\leq p_1+2.$

\textbf{Subcase 2.1.}  $d_G(v_{i+2},y)= p_1+2.$

In this case, $P_2$ is a shortest path from $y$ to $v_{i+2}$. Let
$P_3$ be a second shortest path between $y$ and $v_{i+2}$,
$C_1=P_2\bigtriangleup P_3$, $C_1\cap P_2\cap P_3=\{x',y'\}$. By
Lemma \ref{lem1}, $\pi(v_{i+2},y)\geq 1$ except for the case that
$C_1$ is an odd cycle and $d_{C_1}(x',y')=1$. In this case, the
length of $P_3$ is $(p_1+2)+|C_1|-2=p_1+|C_1|$, which is not less
than $p_1+k$. Consider the path $yP_1v_iv_{i-1}v_{i-2}\cdots
v_{i+2}$. It is a path between $y$ and $v_{i+2}$, and its length is
$p_1+(k-2)<p_1+k$, contrary to the choice of $P_3$.

\textbf{Subcase 2.2.}  $p_1\leq d_G(v_{i+2},y)< p_1+2.$

Let $P_2'$ be a shortest path from $y$ to $v_{i+2}$, and $P_3'$ a
second shortest path between $y$ and $v_{i+2}$. Let
$C_1'=P_2'\bigtriangleup P_3'$, $C_1'\cap P_2'\cap P_3'=\{x',y'\}$.
If $P'_3 = P_2,$ since $g \geq 5$ and $ |E(P'_2)| \geq |E(P_1)|,$
then $d_{C'_1}(x',y') \geq 2,$ and by Lemma \ref{lem1} we have
$\pi(v_{i+2},y)\geq 1$. If $P'_3 \neq P_2$, by Lemma \ref{lem1},
$\pi(v_{i+2},y)\geq 1$ except for the case that $C_1'$ is an odd
cycle and $d_{C'_1}(x',y')=1$. But, this case cannot happen because
the length of $P'_3$ is $|E(P'_2)|+|C_1'|-2\geq p_1+|C_1'|-2\geq
p_1+k-2\geq p_1+3$, which is larger than the length of $P_2$,
contrary to the choice of $P'_3$.

No matter which cases happen, we always have $\pi(v_{i+2},y)\geq 1$.
Similarly, we have $\pi(v_{i-2},y)\geq 1$. Because $k\geq 5$,
$v_{i-2}$ is different from $v_{i+2}$. For all the remaining
vertices in $C$, $\pi(v_j,y)\geq 0$ for $j\neq i-2,i+2$. Then, for a
fixed $y\in V(G)\backslash V(C),$ we get that $\sum_{x\in
V(C)}\pi(x,y)\geq 2$. Therefore,
$$\sum_{x\in V(C), y\in V(G)\backslash V(C)}\pi(x,y)\geq 2(n-k).$$

\textbf{Case 3.}  $x, y\in V(G)\backslash V(C)$.

In this case, $\pi(x,y)\geq 0$.

From above cases, we have
\begin{eqnarray*}
& & Sz(G)-W(G) \\
& = & \displaystyle \sum_{\{x,y\} \subseteq V(G)} \pi(x,y)\\
& =& \displaystyle\sum_{\{x,y\} \subseteq V(C)}\pi(x,y)
+\sum_{\substack{x\in V(C)\\ y\in
V(G)\backslash V(C)}}\pi(x,y)+\displaystyle \sum_{\{x,y\} \subseteq V(G)\backslash V(C)} \pi(x,y) \\
& \geq &{k\choose 2}-k+2(n-k)\\
& =& 2n+\frac{1}{2}k(k-7)\\
& \geq &2n-5.
\end{eqnarray*}
for $k\geq 5.$

From the above inequalities, we see that equality holds if and only
if $k=g=5$, $\pi(x,y)=1$ for all the nonadjacent pairs $\{x,y\}$ in
$C$, and there are exactly two vertices $v_1,v_2$ in $C$ such that
$\pi(v_1,y)=1,\pi(v_2,y)=1$ for all $y \in V(G)\backslash V(C)$, and
$\pi(x,y)=0$ for every pair $\{x,y\}\subseteq V(G)\backslash V(C)$.

We first claim that if the equality holds, then $G$ is unicyclic.
Suppose that $\mathcal{C}$ is the set of all cycles except the
shortest cycle $C$. Let $C'$ is a shortest cycle of $\mathcal{C},$
then $C'$ is an isometric cycle. If $C'$ is an even cycle, and there
exists a pair of vertices $\{x, y\}\subseteq V(C')\backslash V(C)$,
then by Lemma \ref{lem2}, $\pi(x,y)\geq 1$, a contradiction. So
there is only one vertex $x\in V(C')\backslash V(C)$. Let $v_i, v_j$
be the neighbors of $x$ in $C'$. Then $v_ix, xv_j$ together with a
shortest path between them in $C$ is a cycle $C''$ different from
$C$. Since the length of $C$ is $5$, $d(v_i,v_j)\leq 2$, and the
length of $C''$ is at most $4$, contrary to the assumption that
$g\geq 5$.

If $C'$ is an odd cycle, and there exists a pair of nonadjacent
vertices $\{x, y\}\subseteq V(C')\backslash V(C)$. Then by Lemma
\ref{lem2}, $\pi(x,y)\geq 1$, a contradiction. If there are only two
adjacent vertices $x, y$ on $V(C')\backslash V(C)$, and let $v_i$ be
the neighbor of $x$ in $C$ and $v_j$ the neighbor of $y$ in $C$,
then $v_ixyv_j$ together with a shortest path between them in $C$ is
a cycle $C_1$ different from $C$. Since the length of $C$ is $5$ and
$g\geq 5$, $d(v_i,v_j)=2$. Then $C_1$ is an isometric cycle, and by
Lemma \ref{lem2}, $\mu_{v_i,v_j}(xy)=1$, and so $\pi(v_i,v_j)\geq
2$, a contradiction. If there is only one vertex $x \in
V(C')\backslash V(C)$, and let $v_i, v_j$ be the neighbors of $x$ in
$C'$, then $v_ix, xv_j$ together with a shortest path between them
in $C$ is a cycle $C_2$ different from $C$. Since the length of $C$
is $5$, $d(v_i,v_j)\leq 2$, and the length of $C_2$ is at most $4$,
contrary to the assumption that $g\geq 5$.

So, we have that $G$ is a unicyclic graph with the only cycle $C$ of
length $5$. Let $C=v_1v_2\cdots v_5v_1$, $T_i$ be the component of
$E(G)\backslash E(C)$ that contains the vertex $v_i (1\leq i\leq
5)$.

If there are at least three nontrivial $T_i$s, say $T_i, T_j, T_k$,
then there is a pair of vertices, say $\{v_i, v_j\}$ which are not
adjacent. Let $x\in V(T_i)\backslash \{v_i\}$, $y\in
V(T_j)\backslash \{v_j\}$. Then $\{x,y\}\subseteq V(G)\backslash
V(C)$. Since $d_C(v_i,v_j)= 2$, by Lemma \ref{lem1}, $\pi(x,y)\geq
1$, a contradiction. Therefore, there are at most two nontrivial
$T_i$s, say $T_i, T_j$. If $v_i, v_j$ are not adjacent, similarly we
can find $\{x,y\}\subseteq V(G)\backslash V(C)$ satisfying
$\pi(x,y)\geq 1$, a contradiction. Thus, $v_i, v_j$ must be
adjacent. In this case, for any $x\in V(T_i)\backslash \{v_i\}$,
$y\in V(T_j)\backslash \{v_j\}$, $\pi(x,y)=0$, and for any $x\in
V(T_i)\backslash \{v_i\}$, $\pi(x,v_{i-2})=1$, $\pi(x,v_{i+2})=1$,
and $\pi(x,v_k)=0$ for $k\neq i,j$. $y\in V(T_j)\backslash \{v_j\}$
is similar to the $x$ case. By calculation, we have $Sz(G)-W(G) =
2n-5$. If there is only one nontrivial $T_i$, we also can calculate
that $G$ satisfies $Sz(G)-W(G) = 2n-5$.
\begin{qed}
\end{qed}

Here we notice that by the above same way, we can give another proof
to Conjecture \ref{conj1}, and get the following result:
\begin{thm}\label{thm2}
Let $G$ be a bipartite connected graph with $n \geq 4$ vertices and
$m \geq n$ edges. Then
$$
Sz(G)-W(G) \geq 4n-8.
$$
Equality holds if and only if $G$ is composed of a cycle on $4$
vertices $C_4$ and a tree $T$ on $n-3$ vertices sharing a single
vertex.
\end{thm}

\pf Let $C$ be a shortest cycle of $G$, and assume that
$C=v_1v_2\cdots v_gv_1$. Simply, $C$ is an isometric cycle. We
consider the pair $\{x,y\}\subseteq V(G)$.

\textbf{Case 1.}  $\{x,y\}\subseteq V(C)$.

By Lemma \ref{lem2}, $\pi(x,y)\geq d_C(x,y).$ Thus, if $xy$ is an edge of $G$, then $\pi(x,y)\geq 1$. Otherwise, $\pi(x,y)\geq 2$.
Therefore, $$\sum_{\{x,y\}\subseteq V(C)}\pi(x,y)\geq g+2\left({g\choose 2} -g\right).$$

\textbf{Case 2.}  $x\in V(C), y\in V(G)\backslash V(C)$.

Assume that $v_i$ is a vertex on $C$ such that $d_G(v_i,y)=min_{v\in
V(C)}d_G(v,y)$, and $P_1$ is a shortest path between $v_i$ and $y$.
Then $P_1$ does not contain any vertices on $C$; otherwise, if
$v_j\in P_1$, then $d_G(v_j,y)<d_G(v_i,y)$, contrary to the choice
of $v_i$.

If there is only one path between $y$ and $v_i$, then $\pi(y,v_i)=0$
and $v_i$ is a cut vertex. For any other vertex $v_j$ in $C$, the
path from $y$ to $v_j$ must go through $v_i$, and thus,
$\mu_{v_i,v_j}(e)=\mu_{y, v_j}(e)$ for $e\in E(C)$. From Lemma
\ref{lem2}, we have that if $v_iv_j$ is an edge of $C$, then
$\pi(y,v_j)\geq 1$. If $d_C(v_i,v_j)\geq 2$, then $\pi(y,v_j)\geq
2$. Therefore, $$\sum_{x\in V(C)}\pi(x,y)\geq 2+2(g-3)=2g-4\geq g.$$

If there are at least two paths between $y$ and $v_i$, then, since
$G$ is a bipartite graph, by Lemma \ref{lem1} $\pi(y,v_i)\geq 1$.
And for each $v_j\in V(C)\backslash \{v_i\}$, there are at least two
paths from $y$ to $v_j$, so $\pi(y,v_j)\geq 1$. Therefore,
$$\sum_{x\in V(C)}\pi(x,y)\geq  g.$$

\textbf{Case 3.}  $x\in V(G)\backslash V(C), y\in V(G)\backslash V(C)$.

In this case, $\pi(x,y)\geq 0$.

From the above cases, we have
\begin{eqnarray*}
& & Sz(G)-W(G) \\
& = & \displaystyle \sum_{\{x,y\} \subseteq V(G)} \pi(x,y)\\
& =& \displaystyle\sum_{\{x,y\} \subseteq V(C)}\pi(x,y)
+\sum_{\substack{x\in V(C)\\ y\in
V(G)\backslash V(C)}}\pi(x,y)+\displaystyle \sum_{\{x,y\} \subseteq V(G)\backslash V(C)} \pi(x,y) \\
& \geq &g+2({g\choose 2}-g)+g(n-g)\\
& =&g(n-2)\\
& \geq &4n-8.
\end{eqnarray*}

From the above inequalities, one can see that if equality holds,
then $g=4$, and $\pi(x,y)=1$ for all the adjacent pairs $\{x,y\}
\subseteq V(C)$, $\pi(x,y)=2$ for all the nonadjacent pairs $\{x,y\}
\subseteq V(C)$ and $\pi(x,y)=0$ for every pair $\{x,y\}\subseteq
V(G)\backslash V(C)$.

Now we show that if equality holds, then $G$ is a unicyclic graph.
Suppose that $\mathcal{C}$ is the set of all cycles except the
shortest cycle $C$. Let $C'$ is a shortest cycle of $\mathcal{C}.$
Then $C'$ is an isometric cycle. Since $G$ is bipartite, $C'$ is an
even cycle. If there exists a pair of vertices $\{x, y\}\subseteq
V(C')\backslash V(C)$, then by Lemma \ref{lem1}, $\pi(x,y)=1$, a
contradiction. So there is only one vertex $x\in  V(C')\backslash
V(C)$. Let $v_i, v_j$ be the neighbors of $x$ in $C'$. Then $v_ix,
xv_j$ together with a shortest path between them in $C$ is a cycle
$C''$ different from $C$. Since the length of $C$ is $4$,
$d(v_i,v_j)= 2$, and the length of $C''$ is $4$,
$\mu_{v_i,v_j}(xv_i)=\mu_{v_i,v_j}(xv_j)=1$. Thus, $\pi(v_i,v_j)\geq
4$, a contradiction. Therefore, $G$ is unicyclic.

Let $T_i$ be the component of $E(G)\backslash E(C)$ that contains the vertex $v_i (1\leq i\leq 4)$.

If there are at least two nontrivial $T_i$s, say $T_i, T_j$, and let
$x\in V(T_i)\backslash \{v_i\}$, $y\in V(T_j)\backslash \{v_j\}$,
then $\{x,y\}\subseteq V(G)\backslash V(C)$, and there are at least
two paths between $x$ and $y$. By Lemma \ref{lem1}, $\pi(x,y)\geq
1$, a contradiction. Therefore, there is only one nontrivial $T_i$.
In this case, we can calculate that $G$ satisfies $Sz(G)-W(G) =
4n-8$. Hence, equality holds if and only if $G$ is the graph
composed of a cycle on $4$ vertices, $C_4$, and a tree $T$ on $n-3$
vertices sharing a single vertex.

\begin{qed}
\end{qed}

Since for a bipartite graph, we have $Sz^*(G)=Sz(G)$, which
immediately implies Conjecture \ref{conj2}.

Next, we give a proof to Conjecture \ref{conj4}. At first we need
the following Lemmas.

\begin{lem}(\cite{SGB})\label{lem3}
 For a connected graph $G$ with at least two vertices,
 $$
 Sz(G) \geq W(G),
 $$
 with equality if and only if each block of $G$ is a complete graph.
\end{lem}

\begin{lem}\label{lem4}
Let $G$ be a connected graph with $n\geq 4$ vertices and $m\geq n$
edges and with an odd cycle. Then for every vertex $u \in V(G),$
there exists an edge $e=v_1v_2 \in E(G)$ such that $u \in N_0(e)$,
that is, $\sum_{e \in E(G)}n_0(e) \geq n.$
\end{lem}

\pf Suppose that there is a vertex $u \in V(G)$ such that for every
$e=xy \in E(G)$, we have $d_G(u,x)\neq d_G(u,y).$ Let $d=ecc(u)$,
$N^i(u)=\{v\in V(G)|d_G(u,v)=i\}, 1 \leq i\leq d.$ By the
assumption, we know that there is no edge in $N^i(u)$ for every $i$,
that is, $N^i(u)$ is an independent set. Set $X= \{u\} \bigcup
\bigcup _{1 \leq i \leq d, i \mbox{\small \ is even}}N^i(u)$, $Y=
\bigcup_{1 \leq i \leq d, i \mbox{\small \ is odd}}N^i(u)$. Then
$G=G[X,Y]$ is a bipartite graph with partite sets $X$ and $Y$. But,
$G$ is a connected graph with an odd cycle, a contradiction. Hence,
for every vertex $u \in V(G),$ there exists an edge $e=v_1v_2 \in
E(G)$ such that $u \in N_0(e)$, and so we have $\sum_{e \in
E(G)}n_0(e) \geq n.$

\begin{qed}
\end{qed}

Now we turn to solving Conjecture \ref{conj4} and get the following
result:
\begin{thm}\label{thm3}
Let $G$ be a connected nonbipartite graph with $n \geq 4$ vertices.
Then
$$
Sz^*(G)-W(G) \geq \frac{n^2+4n-6}{4}.
$$
Equality holds if and only if $G$ is composed of a cycle on $3$
vertices, $C_3$, and a tree $T$ on $n-3$ vertices sharing a single
vertex.
\end{thm}

\pf By using $n_u(e)+n_v(e)+n_0(e)=n$ for every $e \in E(G)$, we
have
\begin{eqnarray*}
& & Sz^*(G)-W(G) \\
& = & \displaystyle \sum_{e=uv\in E(G)} \left(n_u(e)+\frac{n_0(e)}{2}\right)\left(n_v(e)+\frac{n_0(e)}{2}\right)-W(G)\\
& = & \displaystyle \sum_{e=uv\in E(G)} n_u(e)n_v(e)+\displaystyle \sum_{e=uv\in E(G)}\left(\frac{n_0(e)}{2}(n-n_0(e))+\frac{n^2_0(e)}{4}\right)-W(G)\\
& = & Sz(G)-W(G)+\displaystyle \sum_{e=uv\in E(G)}\left(\frac{n_0(e)}{2}n-\frac{n^2_0(e)}{4}\right)\\
% & \geq & \frac{n^2}{2}-\frac{1}{4}\left(1^2+\cdots+1^2+(n-(g-1))^2\right)\\
% & \geq & \frac{n^2}{2}-\frac{1}{4}\left(2+(n-2)^2\right)\\
% & = & \frac{n^2+4n-6}{4}.
\end{eqnarray*}

Let $n_0=\displaystyle \sum_{e=uv\in
E(G)}\left(\frac{n_0(e)}{2}n-\frac{n^2_0(e)}{4}\right)$. If there
are two edges $e',e''$ such that $n_0(e')\geq n_0(e'')$, and put
$n'_0(e')=n_0(e')+1,$ $n'_0(e'')=n_0(e'')-1,$ $n'_0(e)=n_0(e)$ for
other edges, then
\begin{eqnarray*}
& & n_0'-n_0 \\
& = & \displaystyle \sum_{e=uv\in E(G)}\left(\frac{n'_0(e)}{2}n-\frac{n'^2_0(e)}{4}\right)
-\displaystyle \sum_{e=uv\in E(G)}\left(\frac{n_0(e)}{2}n-\frac{n^2_0(e)}{4}\right)\\
& = & \frac{n_0(e'')-n_0(e')-1}{2}\\
& < & 0.
\end{eqnarray*}

Let $C$ be a shortest odd cycle of $G$ with length $g$, and its
edges be $e_1, e_2,\cdots, e_g$. Then $C$ is isometric. For every
edge $e=uv\in E(C)$, there is a vertex $x\in V(C)$ such that
$d_G(x,u)=d_C(x,u)=d_C(x,v)=d_G(x,v)$. Therefore, $n_0(e)\geq 1$ for
every $e\in E(C)$. If there are two edges $e',e''$ such that
$n_0(e')\geq n_0(e'')$, we could do the operation as above, which
makes $n_0$ smaller. Thus, $n_0$ attains its minimum when
$n_0(e_i)=1$ except for $n_0(e_1)$, $n_0(e)=0$ for all the remaining
edges. By Lemma \ref{lem4}, $\sum_{e \in E(G)}n_0(e) \geq n$, and so
$n_0(e_1)\geq n-g+1$. Hence, $$n_0\geq
(g-1)(\frac{n}{2}-\frac{1}{4})+\frac{n-g+1}{2}n-\frac{(n-g+1)^2}{4}
 \geq \frac{n^2}{2}-\frac{1}{4}\left(2+(n-2)^2\right)
=\frac{n^2+4n-6}{4}. $$

From the above inequalities, we can see that equality holds if and
only if $g=3$, $Sz(G)=W(G)$ and $n_0(e_1)=n-2, n_0(e_2)=1,
n_0(e_3)=1, n_0(e)=0$ for all the remaining edges.

Now we conclude that $G$ is unicyclic. Suppose that $G$ is not
unicyclic. By Lemma \ref{lem3}, we know there is a block $H$
different from $C$ which is a complete graph of order at least
three. Then, $n_0(e) \geq 1$ for every $e \in E(H)$, a
contradiction.

Let $T_i$ be the component of $E(G)\backslash E(C)$ that contains the vertex $v_i (1\leq i\leq 3)$.

If there are at least two nontrivial $T_i$s, say $T_1, T_2$, then
$n_0(v_2v_3)= |V(T_1)| \geq 2, n_0(v_1v_3)= |V(T_2)| \geq 2$, a
contradiction. Therefore, there is only one nontrivial $T_i$. In
this case, we can calculate that $G$ satisfies $Sz^*(G)-W(G) =
\frac{n^2+4n-6}{4}$. Hence, equality holds if and only if $G$ is the
graph composed of a cycle on $3$ vertices, $C_3$, and a tree $T$ on
$n-3$ vertices sharing a single vertex.
\begin{qed}
\end{qed}

\end{document}